\RequirePackage[british]{babel}
\documentclass[a4paper]{amsart}
\usepackage[T1]{fontenc}
\usepackage{amscd}
\usepackage{url}

\theoremstyle{plain}
\newtheorem{thm}{Theorem}[section]
\newtheorem{cor}[thm]{Corollary}
\newtheorem{prop}[thm]{Proposition}
\newtheorem{lem}[thm]{Lemma}

\theoremstyle{definition}

\newtheorem{ex}{Example}

\theoremstyle{remark}
\newtheorem*{rmk}{Remark}
\newtheorem*{acks}{Acknowledgments}

\newcommand{\ie}{\textit{i.e. }}

\newcommand{\p}{\mathbb{P}}

\newcommand{\T}{\mathbb{T}}

\newcommand{\Gr}{\mathbb{G}}

\newcommand{\abs}[1]{\lvert#1\rvert}

\newcommand{\gen}[1]{\langle#1\rangle}

\DeclareMathOperator{\length}{length}

\DeclareMathOperator{\red}{red}

\DeclareMathOperator{\Sing}{Sing}
\DeclareMathOperator{\im}{Im}
\DeclareMathOperator{\jo}{J}

\begin{document}

\title[Congruences in $\p^4$ with Non-reduced Fundamental Surface]{On First 
Order Congruences of Lines in $\p^4$ with Generically 
Non-reduced Fundamental Surface}

\author{Pietro De Poi}
\thanks{This research was partially supported by the DFG Forschungsschwerpunkt
``Globalen Methoden in der Komplexen Geometrie'', and the EU (under the EAGER
network).}
\address{Dipartimento di Matematica e Informatica\\
Universit\`a degli studi di Trieste\\
Via Valerio, 12/1\\
34127 Trieste\\
Italy\\
}
\email{depoi@dmi.units.it}

\date{\today}

\begin{abstract}
In this article we obtain a complete description of the  congruences of lines 
in $\p^4$ 
of order one provided that the fundamental surface $F$ is non-reduced 
(and possibly reducible) at one of its generic points, and 
their classification under the hypothesis that $(F)_{\red}$ is smooth. 
\end{abstract}

\maketitle

\bibliographystyle{amsalpha}

\section*{Introduction}
A \emph{congruence} of $d$-dimensional linear subspaces in $\p^n$ is 
an irreducible subvariety (``family'') $B'$ of dimension $n-d$ of the 
projective Grassmannian $\Gr(d,n)$. By taking a resolution $B$ of 
singularities of $B'$, one can pull back the universal subbundle on 
$\Gr(d,n)$, obtaining a smooth $\p^d$-bundle $\Lambda$ on $B$, which has 
dimension $n$. Since $\Lambda\subset B\times \p^n$, the second projection 
induces a morphism $p\colon\Lambda\rightarrow \p^n$ between manifolds of the 
same dimension, and one defines the \emph{order} of the family as the 
degree of $p$. Geometrically, this is the number of $d$-planes of the family 
passing through a general point in $\p^n$. 

The case where $B$ is a surface ($n-d=2$), and the order is one, was classified 
by Ziv Ran (\cite{R}), extending the classical work of Kummer (\cite{K}) 
for $n=3$ and order at most two. 

The case where $n=4$, $d=1$ and the order is one was considered by G. Marletta 
(\cite{M1}, \cite{M3}) and, since  \cite{tesi}, we are trying to complete 
Marletta's incomplete classification and to bring it up to modern standards. 
In fact, in \cite{tesi} we have started the 
classification of first order congruences of lines in $\p^4$, 
and the first two steps of the classification are in 
in \cite{DPi} and \cite{DP2}.

In this article, we study one of the most difficult and exciting cases, 
namely when the fundamental surface (see below) of our congruence is non-reduced 
at one of its generic points. 
The study of these cases is not only interesting in its own and for its 
applications to projective geometry, but it will
also be useful---hopefully---for mathematical physics, in the theory of  
systems of conservation laws, see  
\cite{AF} and \cite{AF1}. In particular, the results contained here can be 
used to classify the three dimensional 
Temple systems of conservation laws which are hyperbolic 
but \emph{not} strictly hyperbolic, and so, among other 
things, can complete the 
results contained in \cite{AF2}. 
 
From now on, we will consider always, for simplicity, $d=1$. 
The set up for studying congruences of lines is the following: let us denote 
by $R$ the ramification divisor of $p$, and by $\Phi$ its schematic image, 
\ie the branch locus, which will be called the \emph{focal locus}. 

We can observe that, if the order is at least two, by the purity of 
the branch locus, $\Phi$ contains a hypersurface, whereas, if the order is 
one, the morphism $p$ is birational. Thus, by Zariski's main theorem, each  
component of the focal locus has codimension two, 
since it is precisely the locus where 
the fibre dimension is positive; moreover, set-theoretically, $\Phi$  
coincides with the \emph{fundamental locus}, \ie the set of points in 
$\p^n$ through which there pass infinitely many lines of the congruence.

The \emph{pure fundamental locus} $(\emptyset\neq)F\subset\Phi$ is instead 
the image under $p$ of the 
components of $R$ which surjectively project onto $B$; this means geometrically 
that every line of $B$ intersects each component of $F$.

If the order is one and $n=4$, 
it can be proven that either $\dim(\Phi)=0$, in which 
case $B$ is a star of lines 
(\ie the set of lines passing through the point $(\Phi)_{\red}$) or 
$\dim(\Phi)=2$ (see Theorem 2.2 in \cite{DP4}). 
In what follows, we shall only consider congruences of lines in 
$\p^4$ such that $\Phi$ has 
pure dimension two. 
In this situation  
$B'$ is a three-dimensional 
subvariety of the trisecant lines to $F$, by a result of C. Segre  
(see Proposition~\ref{prop:fofi} below). 
The cases in which $F$ contains a component of 
dimension one is treated 
in  \cite{DP2}, while the case in which the pure fundamental locus $F$ is a 
(generically) reduced and irreducible surface is in \cite{DPi}. 

Here we are interested in the cases in which $F$ has some component of 
dimension two such that $F$ is nonreduced at the corresponding generic 
point. $F$ will also be called \emph{fundamental surface} and it is  
either irreducible or reducible. 
For both these possibilities we give a precise description.

Classically, congruences of order one in $\p^4$ were treated by G. Marletta 
in \cite{M1} 
for the generically reduced case and in \cite{M3} for the rest. 
In this paper 
Marletta's classifications are reproved and completed.  
In particular 
we deduce the complete 
list of these congruences provided that 
all the reduced components of $F$ are smooth. 

This article is structured as follows: 
after giving, in Section~\ref{sec:1}, the 
basic definitions, we give examples of all the possible congruences with a non-reduced component 
of the fundamental locus  (finding also a case 
which is missing in Marletta's list) in Section~\ref{sec:2b}, and in 
Sections~\ref{sec:2} and \ref{sec:3} 
we prove that no other congruence with this characteristic property can exist. 
A partial list of these congruences has also been 
obtained by A.~Oblomkov (\cite{O}). If we suppose that all the components of the 
reduced locus of the pure fundamental locus are smooth, 
we obtain the complete list in Theorem~\ref{thm:1er}. 

In this theorem and throughout  
this article we will use the following notation and conventions: 
$F$ will be the pure fundamental locus, and moreover  
we set $D:=(F_1)_{\red}$, where $F_1$ is a two-dimensional 
component of $F$ such that $F$ is nonreduced at the generic point of $F_1$. 
If $b\in B$, then $\Lambda(b)\subset\p^4$ is the corresponding line of the 
congruence.

\begin{thm}\label{thm:1er}
If the fundamental surface $F$ of a first order congruence of lines in  
$\p^4$ is such that $F_1\neq \emptyset$  
and all the components of $(F)_{\red}$ 
are smooth, then we have the following possibilities: 
\begin{enumerate}

\item $F(=F_1)$ is irreducible, and we have the following cases:

\begin{enumerate}
\item $\length(\Lambda(b)\cap D)=1$, 
in which case $D$ is a plane and the 
congruence is  as in Example~\textup{\ref{ex:1}}; or

\item\label{1:1}  $\length(\Lambda(b)\cap D)=2$,  $D$ 
is a rational normal cubic scroll, and the congruence is 
as in Example~\textup{\ref{ex:2}}; or
\end{enumerate}

\item\label{1er} $F$ has two irreducible components, $F_1$ non-reduced 
and $F_2$ reduced, and we have the following cases:

\begin{enumerate}
\item $D=(F_1)_{\red}$ is a plane and $F_2$ is a rational normal cubic 
scroll, and we have the following cases:
\begin{enumerate}
\item $D\cap F_2$ is either a line or a conic and the congruence is 
as in Example~\textup{\ref{ex:4}}; or

\item $D\cap F_2$ is a (smooth) conic  and the congruence is as in Example~\textup{\ref{ex:3}}; or
\end{enumerate}

\item $F_2$ is a plane, and $D=(F_1)_{\red}$ is a rational normal cubic 
scroll; $D\cap F_2$ is a smooth conic, unisecant to $D$, and the congruence is 
as in Example~\textup{\ref{ex:5}}.  
\end{enumerate}
\end{enumerate}

\textit{Vice versa} a family of lines of dimension three 
constructed as in each of the cases of the 
theorem is a first order congruence in $\p^4$.
\end{thm}

\begin{acks}
This article has been highly improved by many discussions with and suggestions 
of F. Catanese.
I want to thank him, E. Mezzetti and E. Arrondo
for their help in the elaboration of this paper.  

I would also like to thank A.~A.~Oblomkov for the useful discussions on the 
subject.
\end{acks}

\section{Notation, Definitions and Preliminary Results}\label{sec:1}
We will work with schemes and varieties over the complex field $\mathbb C$, with 
standard notation and definitions as in \cite{H}. In this article, a 
\emph{variety} will always be projective and irreducible. 
We refer to \cite{DP2} and \cite{DP3} for general results and references 
about families of lines, focal diagrams and congruences, and to \cite{GH} for 
notations about Schubert cycles. So, we denote  by $\sigma_{a_0,a_1}$ 
the Schubert cycle of the lines in $\p^4$ contained in a fixed  
$(4-a_1)$-dimensional subspace $H\subset\p^4$ and which meet a fixed $(3-a_0)$-dimensional 
subspace $\Pi\subset H$. 
We recall that
a \emph{congruence of lines} in $\p^4$ is a (flat) family 
$(\Lambda,B,p)$ of lines 
in $\p^4$ obtained as the pull-back of the universal family under the desingularization of
a subvariety $B'$ of dimension 
three of the Grassmannian $\Gr(1,4)$. 
So, $\Lambda\subset B\times \p^4$ and $p$ is the restriction of the projection
$p_1\colon B\times\p^4\rightarrow B$ to
$\Lambda$, while we will denote the restriction 
of $p_2\colon B\times\p^4\rightarrow\p^4$ by $f$. $\Lambda_b:=p^{-1}(b)$, $(b\in B)$ is 
a line of the family and $f(\Lambda_b)=:\Lambda(b)$ is a line in $\p^4$. 
$\Lambda$ is smooth of dimension three: therefore we can define the \emph{focal divisor} $R\subset \Lambda$ 
as the ramification divisor of $f$.      
The \emph{focal locus}
$\Phi=f(R)\subset \p^4$, is the scheme theoretic image of $R$. 

In this article, we study the \emph{first order congruences} of lines,  
\ie congruences $B$ for which through a general point in $\p^4$ there passes 
only one line of $B$; or, equivalently, we can write  
\begin{equation}\label{eq:class}
[B']=\sigma_3+a\sigma_{2,1}, 
\end{equation} 
where $[B']$ is the rational equivalence class of $B'$---as a linear combination of 
Schubert cycles of the Grassmannian; so, the \emph{class} $a$  
is the degree of the ruled surface 
generated by the lines of the congruence which belong
to a general $\p^3$ (\ie the intersection number 
$[B']\cdot \sigma_{2,1}$). 
Given a first order congruence $B$, through a focal 
point there will pass infinitely many lines in $B$, \ie a focal point 
is a \emph{fundamental point} and the focal locus coincides set-theoretically 
with the \emph{fundamental locus}.
 
An important result, due to Corrado Segre,  is the following (see \cite{DP2}, Proposition~1 for a proof):
\begin{prop}[C. Segre, \cite{S}]\label{prop:fofi}
On every line $\Lambda(b)\subset\p^4$ of the family, 
the  focal locus $\Phi$
either coincides with the whole $\Lambda(b)$---in which case $\Lambda(b)$ 
is called \emph{focal line}---or is a zero dimensional 
scheme of $\Lambda(b)$ of length three. 
\end{prop}

Let us observe that, if the order is one, the morphism $p$ is birational, thus, 
by Zariski main theorem (finite plus birational cover of a normal variety is 
an isomorphism), the focal locus is not a hypersurface, so $\dim(\Phi)\le 2$. 
Actually, we have proven in \cite{DP4} 
that $\dim(\Phi)=2$ if $B$ is not a star of lines. 

Here we are interested in first order congruences of lines in $\p^4$  for 
which the focal locus $\Phi$ has pure dimension two, so 
the ``important'' component of $\Phi$ is the \emph{fundamental $2$-locus} 
(see \cite{DP2}, \cite{DP3} and \cite{DP4}), which is characterized 
by the fact that it is formed by the components of $\Phi$ of pure dimension 
two  such that the general line of the congruence meets it in a zero 
dimensional scheme (see \cite{DP2}, 
Proposition~2). 
The fundamental $2$-locus will be called 
\emph{pure fundamental locus} or, in what follows, 
simply \emph{fundamental surface} 
and it will be indicated with $F$ 
(For the other ``non-important'' components of $\Phi$, \ie the \emph{parasitic
planes}, see \cite{DPi}). 

In the rest of the paper, we need more notation. To a first order congruence $B$ we associate: 
\begin{itemize}
\item the hypersurface $V_\Pi$ in 
$\p^4$ given by the lines of $B$ which intersect a general plane 
$\Pi\subset\p^4$. So, $\deg(V_\Pi)=1+a$, where $a$ is the class of $B$ (see \eqref{eq:class});   
\item the surface $\Sigma_H$---of degree $a$---of the lines of  
$B$ contained in a general hyperplane $H$; We can think of it as the image 
of the $\p^1$-bundle over the curve $\Gamma_H\subset B$, obtained by pulling 
back the intersection of $B'$ with the Schubert variety of the lines contained 
in $H$. 
\item we set also $C_H:=(F)_{\red}\cap H(\subset \Sigma_H)$. 
\end{itemize}

\begin{rmk}
First of all we note that $\Sing(\Sigma_H)\subset C_H$, since 
through a singular point $P$ of the surface there pass more than a line of the 
congruence (possibly infinitely near). Besides, for the generality of $H$ and 
dimensional reasons, $\Sigma_H$ cannot be contained in $F$.
\end{rmk}

\begin{lem}\label{lem:cat}
Let $D$ be a surface in $\p^4$ and $P$ a smooth point in $D$ such that,
for a general point $Q$ in $D$, the tangent planes $\T_{P,D}$ and $\T_{Q,D}$
do not span  $\p^4$. Then the surface $D$ is degenerate. 
\end{lem}

\begin{proof} Take projective coordinates such that the local parametrization of
$D$
at $P$ is $F:= (1:x:y:f(x,y):g(x,y))$, where $f,g$ vanish of order $\geq 2$
for $x=y=0$. Our statement is that the 6 vectors $e_0, e_1, e_2$ and
$F, F_x, F_y$ have rank (at most) 4. This simply means that the three vectors
$(f,g), (f_x,g_x), (f_y,g_y)$ have rank (at most) 1, and implies that the map
$(x:y) \mapsto (f:g)$ has constant image in $\p^1$. Without loss of
generality we may then assume $ g \equiv 0$, hence $D$ is degenerate.
\end{proof}

\begin{prop}
The congruence $B$ of the tangent lines to a surface $D\subset\p^4$ cannot
have order one.
\end{prop}

\begin{proof}
First we note that if $D$ is a cone, $B$ has order zero. Moreover, $D$ is
contained in the focal locus $F$ of $B$. If 
$P\in \T_{P_1,D}\cap \T_{P_2,D}$, then $P\in F$, 
since two lines pass through 
it.

We cannot have $\dim(\T_{P_1,D}\cap \T_{P_2,D})=1$ by the preceding lemma, so 
%
%
we suppose that $\T_{P_1,D}\cap \T_{P_2,D}=:Q$ 
is a point and that $B$
is a first order congruence. 
Now, if $\dim(\T_{P,D}\cap F)=0$ 
for the general $P\in D$, then $Q\in(\T_{P,D}\cap F)$, since  
$\T_{P,D}\cap \T_{P',D}\in \T_{P,D}\cap F$ as $P'$ varies in $D$. 
Therefore $Q$ is 
fixed; but then the dual variety $D^*$ would be degenerate, and $D$ would be a
cone (this argument holds only on a field of characteristic zero). 

Instead $\dim(\T_{P,D}\cap F)=1$ means that there is a 
family of dimension two of planes curves contained in $F$. But the surfaces
with this property 
in $\p^4$ are classified in C.~Segre's Theorem 
(see \cite{MP}, Theorem~4 for a modern proof) and they are 
the projected Veronese surface, the rational normal smooth cubic 
scroll and the cones. 
If $D$ is this surface, we conclude by observing that couples of tangent 
planes to 
either the Veronese surface or the cubic ruled surface do not meet in a 
point in $D$.  

Otherwise, there is another nondegenerate (if it is degenerate, also $D$ is 
degenerate) 
surface $D'\subset F$ which is in the list of  C.~Segre's Theorem, and the 
congruence $B$ is given by the  lines which are tangent to $B$ and these lines moreover
meet  $D'$. In 
particular, on every line $r$ of $B$, the three foci of $r$ are a double point 
in $D$ and a single point in $D'$. $D'$ cannot be a  cone, since again 
$D^*$ would be degenerate; in fact all the tangent planes would pass through 
the vertex of $D'$. Then $D'$ is either the projected Veronese surface or the 
cubic scroll, and in both these cases the plane curves are irreducible 
conics $C\subset \T_{P,D}$. We recall that the lines of $B$ contained in  
$\T_{P,D}$
form the pencil of lines through the point $P\in D$, $\p^1_P$. 
But if $C\subset\T_{P,D}$ is the conic, and 
$\ell_P\in\p^1_P$, then either the two points $\ell_P\cap C$ 
are focal, but then $\ell_P\subset F$ and therefore every line in $B$ is focal, 
or only one of the 
two points is focal for $\ell_P$. But then, we would obtain a regular map 
from $\p^1_P\cong \p^1$ to $C$ which is 
injective not but surjective. 
\end{proof}

\begin{cor}\label{lem:curva}
$\Sigma_H$ does not contain a component $S_H$ 
which is the surface of the tangents to a curve $C$.  
\end{cor}

\begin{proof}
If it were so, $C\subset C_H$ and then 
$B$ would be given by the  
family of tangent lines to $D=(F_1)_{\red}$. 
\end{proof}

\begin{lem}\label{lem:cono}
$\Sigma_H$ does not contain a component $S_H$ 
which is a cone.  
\end{lem}

\begin{proof}
If every $S_H$ would be a cone, since this cone spans the hyperplane $H$, we
would have $\infty^4$ such cones, whence $\infty^5$ pairs $(S_H, L)$, 
where $L$ is a line of the cone $S_H$. Since such a line $L$ belongs to 
$\infty^2$ hyperplanes, it follows
that we get in this way $\infty^3$ lines, hence the whole congruence of
lines.

Consider the irreducible variety $S$ image of the rational map
$H \to V_H$. We know that $\dim(S) \leq 2$, hence a general 
$V = V_H$ is a vertex of at least $\infty^2$ cones, 
and \textit{a fortiori} there are at least $\infty^2$ lines of the
congruence which pass through $V$. Hence we obtain that $S$ is a $1$-dimensional
component of the focal locus, which is excluded by our assumption.


\end{proof}

 
\begin{lem}\label{lem:quadrica}
$\Sigma_H$ does not contain a component $S_H$ 
which is a quadric. 
\end{lem}

\begin{proof}

If $S_H$ is a quadric, it must be smooth, by the preceding lemma. 

If $S_H$ is a smooth quadric, only one of the two families of lines 
contained in it can be contained in $B$, otherwise $S_H\subset F$, 
and this is absurd by varying $H$.

Let $\ell$ be a line of the congruence contained 
in $S_H$. If $\ell'$ is another general line of the congruence not 
contained in 
$H$, the hyperplane $H':=\overline{\ell'\ell}$ determines another quadric 
$S_{H'}$. By construction, $S_H\cap S_{H'}=\ell\cup\ell''$, and 
$\ell'\cap\ell''\neq\emptyset$; therefore $\ell'$ intersects $S_H$, which 
is absurd since every line of $B$ should intersect $S_H$. 
This again implies  $S_H \subset F$, the same
contradiction.
\end{proof}

\section{The examples}
\label{sec:2b}

Since the congruences we are studying are rather complicated to construct, we first give 
some examples describing them, and afterwards we will show that these examples exhaust all 
the possible congruences. 

\begin{rmk}
We observe first that if a component of the focal locus is contained 
set-theoretically in a hyperplane $H$, then 
the lines of the congruence contained in $H$ form a first order congruence. Therefore a way of 
constructing congruences in $\p^4$ (or in general in $\p^n$) is to reverse this by considering 
a one-dimensional linear system of hyperplanes in $\p^4$ such that in each of them there is a first order congruence 
of lines in $\p^3$, and such that these congruences vary algebraically. 
\end{rmk}

Then, we start recalling the classification of first order congruences in $\p^3$: 
we follow the notation used in \cite{DP4}: let $\ell$  be a fixed line in $\p^3$, then $\p^1_\ell$ 
is the set of the planes containing $\ell$. Let $\phi$ be
a general nonconstant morphism from $\p^1_\ell$ to $\ell$ and let $\Pi$ be a general
element in $\p^1_\ell$. 
We define $\p^1_{\phi(\Pi),\Pi}$ as the pencil of lines passing through the point $\phi(\Pi)$ and
contained in $\Pi$.

If $\mathcal{B}$ is a first order congruence of lines in $\p^3$ with focal locus $\mathcal{F}$, then  
\begin{enumerate}
\item $P:=\mathcal{F}_{\red}$ is a point and $\mathcal{B}$ is the star of lines passing through $P$; or 
\item\label{cubica} $\mathcal{F}$ is a rational normal curve in $\p^3$ and $\mathcal{B}$
is the family of secant lines of $\mathcal{F}$; or 
\item\label{retta} $\ell:=\mathcal{F}_{\red}$ is a line,
and the congruence is $\cup_{\Pi\in \p^1_\ell}\p^1_{\phi(\Pi),\Pi}$; or
\item\label{ridu}  $\mathcal{F}=\mathcal{F}_1\cup\mathcal{F}_2$ where
$\mathcal{F}_1$ is a line  and $\mathcal{F}_2$ is a rational curve, such that
$\length(\mathcal{F}_1\cap \mathcal{F}_2)=\deg(\mathcal{F}_2)-1$ and 
$\mathcal{B}$ is the family of lines meeting $\mathcal{F}_1$ and $\mathcal{F}_2$. 
\end{enumerate}

In order to state the examples of congruences in $\p^4$, we need some more notations: 

\begin{itemize}
\item if $H\subset\p^4$ is a linear space, $H^*$ denotes its dual, 
\ie the linear space of the hyperplanes 
in $H$; if $K\subset H$ is another subspace, we set $K^\vee_H:=\{\Pi\in H^*\mid \Pi\supset K \}\subset H^*$ 
(or simply $K^\vee$ if $H=\p^4$) \ie $K^\vee_H$ is the projective dual subspace 
of $K$ in $H$; $k$ will denote an element in $K^\vee_H$; 

\item $\abs{C}$ denotes 
the complete linear system associated to the divisor $C$;

\item a general (rational) map from $X$ to $Y$  will be denoted by $\phi_{X,Y}$; 

\item if $X,Y\subset\p^4$, $\jo(X,Y)\subset\Gr(1,4)$ is 
$\jo(X,Y):=\{\ell\in\Gr(1,4)\mid \ell\cap X\neq\emptyset\neq\ell\cap Y\}$;  
$\overline{X}$ and $\gen{X}$ are, 
respectively, the (Zariski) closure and the span of $X$. 

\item We recall that we set $D:=(F_1)_{\red}$, where $F_1$ is the component of $F$ which is non-reduced. 

\item Finally, in the case in which $F=F_1\cap F_2$ is also reducible, we denote $C_{H,1}:=(F_1)_{\red}\cap H$ and 
$C_{H,2}:=F_2\cap H$, so $C_H=C_{H,1}\cup C_{H,2}$. $k_i$ is the algebraic 
multiplicity of $C_{H,i}$ in $\Sigma_H$. 

\end{itemize}

\begin{ex}\label{ex:1}
Let us see how to construct a congruence of lines in $\p^4$ from the case~\eqref{retta} in $\p^3$, with set-theoretically linear focal locus, \ie $D:=F_{\red}$ 
is a plane.  
We want to follow the remark above, so we need a pencil of hyperplanes, and 
therefore we have to fix a 
plane $D\subset\p^4$, so that the pencil is $D^\vee$. 
Now, for every element $d\in D^\vee$, we have to find a line $\ell_d\subset D$, in order to construct 
the congruence of lines in $d\cong\p^3$. Since we do not want (embedded) components of dimension one, 
we suppose that the line $\ell_d$ is not the same for all the elements in $D^\vee$. 
In brief, we have fixed a morphism $\phi:=\phi_{D^\vee,D^*}$.  

Now, we have to construct a congruence in every
$d\in D^\vee$, as in  case~\eqref{retta}, and therefore we have to fix a morphism  
$\psi_d:=\phi_{\phi(d)^\vee_d,\phi(d)}$; moreover, we are in the algebraic category, and so 
we assume that we have an algebraic family of morphisms 
$\{\psi_d\}_{ d\in D^\vee}$. 

Now it is clear that our congruence in $\p^4$ is formed
by the lines of the pencils  $(\psi_d(h))^\vee_h$ as $h$ varies in
$\phi(d)^\vee_d$  and
$d$ in $ D^\vee$, \ie 
$$
B=\overline{\cup_{d\in D\vee}
\cup_{h\in \phi(d)^\vee_d}(\psi_d(h))^\vee_h};
$$
an easy calculation (or one can see Theorem~8, case~3 of \cite{DP2} where this congruence is 
also introduced) shows that the bidegree of $\Lambda$ is $(1,d_1d_2m+1)$, where
$d_1:=\deg(\phi)$,
$d_2:=\deg(\psi_d)$ and $m:=\deg(\phi(D^\vee))$ and
if $r$ is the general line in $D$, then $\deg(f^{-1}(r))=d_1d_2m$.

We note that we could construct more first order congruences if we associate to 
$d$ a general line in $\p^4$ (and not contained in $D$). Actually, in this case 
the reduced focal locus $D:=F_{\red}$ is a ruled surface 
and it is easy to show that $B$ is given 
by a family of secant lines to $D$, and in fact this case is 
contemplated in the example which follows. 
\end{ex}

\begin{ex}\label{ex:2}
We will see now how to construct a first order congruence of lines $B$ such that it is  
a subfamily of the secant lines 
to the reduced locus $D:=(F)_{\red}$.  Let us suppose that $D$ is 
a rational normal cubic scroll in $\p^4$, \ie either a non-degenerate 
cone or a rational normal (smooth) scroll of type $(1,2)$. 
In order to do so, let us consider a curve $C_0$ which is either 
the zero section (see Proposition~\textup{V.2.8} of \cite{H}) if the scroll is smooth, 
\ie an irreducible conic, or a pair of lines if $D$ is a cubic cone, and then 
let us fix a map 
$\phi_{\abs{C_0},F}$ 
which 
associates to the curve $C\in \abs{C_0}$ a point $P_C:=\phi_{\abs{C_0},F}(C)\in C$; 
then we define a congruence as
$$
B:=\overline{\cup_{C\in \abs{C_0}}(P_C)^\vee_{\gen{C}}}; 
$$
it is clearly a first order congruence of lines and 
moreover, if $k=\deg(\phi_{\abs{C_0},F})$ then the class of $B$ is $2k$. 
\end{ex}

If we want to obtain congruences of lines in $\p^4$ from case~\eqref{ridu}, 
we have in fact many possibilities, 
since we may suppose either that the line or the rational curve is in the non-reduced locus. In particular, we will see 
that in Example~\ref{ex:3} we have the case of the curve (a conic, indeed), 
and in Example~\ref{ex:4} the line; 
in both cases the non-reduced locus is a plane. 

\begin{ex}\label{ex:3}
So, let us first suppose that $D:=(F_1)_{\red}$ is a plane, and that the rational curves of case~\eqref{ridu} are 
conics contained in $D$. Then, we also suppose that 
$F_2$ is a non-degenerate rational scroll of type $(1,d)$ and by what we are supposing, the ruling lines must be the lines of 
case~\eqref{ridu}; therefore, 
$D\cap F_2$ is a unisecant curve (of degree $d$) in $F_2$.  
Now, let $\abs{H_1}\cong \p^1$ be the linear system of 
the lines of the ruling in $F_2$ 
and $\abs{C_2}$ the 
linear system of conics in $D$; 
then the congruence is defined as follows: fix a morphism 
$\phi:=\phi_{\abs{H_1}, \abs(C_2)}$ 
such that $\phi_{\abs{H_1}, \abs{C_2}}(\ell)\cap\ell\neq \emptyset$; then we set  
$$
B:=\overline{\cup_{\ell\in \abs{H_1}}\cup_{r\in \jo(\ell,\phi_{\abs{H_1}, \abs{C_2}}(\ell))}r}; 
$$
we also have that $k_1=(d+1)\deg(\phi)\deg(\im(\phi))$, $k_2=2$, and so 
$a=2+\deg(F_2)\deg(\phi)\deg(\im(\phi))$. 
\end{ex}

\begin{ex}\label{ex:4}
Again, let us first suppose that $D:=(F_1)_{\red}$ is a plane, but now we suppose that in it there are contained the lines 
of case~\eqref{ridu}. Now, we suppose also that 
 $F_2$ is a non-degenerate rational surface with sectional genus 
zero (\ie either a rational scroll or the projected Veronese surface); 
let $d\in D^\vee$ and let $\abs{H_1}$ be the corresponding 
linear system of the curves $H_1:=\overline{(F_2\cap d)\setminus D}$; 
then the congruence is given in this way: fix a morphism 
$\phi:=\phi_{\abs{H_1}, D^*}$ 
such that $\length(\phi_{\abs{H_1}, D^*}(C_2)\cap C_2)= \deg(C_2)-1$; then we set  
$$
B:=\overline{\cup_{C_2\in \abs{H_1}}\cup_{r\in \jo(C_2,\phi_{\abs{H_1}, D^*}(C_2))}r}; 
$$
we have also that $k_1=\deg(F_2)\deg(\phi)\deg(\im(\phi))$ and $k_2=1$, and so 
$a=1+\deg(F_2)\deg(\phi)\deg(\im(\phi))$. 
\end{ex}

In the last example the focal locus is reducible, and the non-reduced 
component will be a rational ruled surface. Also this, as Example~\ref{ex:1}, 
is obtained from case~\eqref{retta} in $\p^3$:

\begin{ex}\label{ex:5}
Now $F_2$ is a plane, 
suppose that $D=(F_1)_{\red}$ is a non-degenerate rational 
scroll of type $(1,c)$ and that  
$D\cap F_2$ is a unisecant curve in $D$ of degree $c$. 
Let $\abs{H_1}\cong \p^1$ 
be the linear system of the lines of the ruling in $D$; 
the congruence is obtained in the following way: 
$\forall L \in \abs{H_1}$,
fix a morphism $\psi:=\phi_{L^\vee_{\gen{L,F_2}},L}$; 
then we set 
$$
B:=\overline{\cup_{L\in  \abs{H_1}} \cup_{l\in L^\vee_{\gen{L,F_2}}}\psi(l)^\vee_l};
$$   
we have also that $k_1=\deg(\psi)$ and $k_2=c+1$, and so 
$a=\deg(D)+\deg(\psi)$. 
\end{ex}

\section{The Case of  an Irreducible Fundamental Surface}
\label{sec:2}
We start considering the case in which $F$ is irreducible; 
since we are interested 
in the cases in which $F$ is non-reduced at a generic point, 
we have two possibilities: 
either the general line $\Lambda(b)$ intersects $F$ in a fat point of 
length three, or it intersects $F$ in two points, one of 
which is a fat point of length two and the other is a simple point.

\subsection{The Case in which $\length(\Lambda(b)\cap (F)_{\red})=1$} 

If the general line of the congruence intersects $F$ in only a fat point,
we can prove the following:

\begin{thm}\label{prop:1}
With notation as in the introduction, 
if  $\length(\Lambda(b)\cap D)=1$, then $D$ is a plane and the 
congruence is given in as in Example~\textup{\ref{ex:1}}. 

\textit{Vice versa} a family of lines constructed in this way is a first order 
congruence. 
\end{thm}

\begin{proof}
Let us consider the surface $\Sigma_H$ of degree $a$ associated to a general 
hyperplane $H$; it contains the curve $C_H$ with some algebraic multiplicity, say $k$. 

If $L\subset\Sigma_H$ is a (general) line, and $l\in L^\vee_H$ a general plane, then 
$l\cap\Sigma_H=L\cup C$,  where $C$ is a curve of degree $a-1$, with a point of multiplicity $k-1$ 
in $L\cap F$. 
$C$ intersects $L$ in only another point, \ie the sole 
point in which $l$ is tangent to the irreducible component $S_L$ of $\Sigma_H$ 
which contains $L$. We recall in fact that 
a general tangent plane to a surface $S$ is tangent at only one point in $S$ 
but in the case in which $S$ is a surface of tangents to a curve or a cone 
(see \cite{Zl}), 
and $S_L$ can be neither a surface of tangents to a curve nor a cone, by 
Corollary~\ref{lem:curva}, and Lemmas~\ref{lem:cono} and \ref{lem:quadrica}; 
therefore we have that $a=k+1$. 

I claim now that 
$C_H$ is a plane curve. Indeed, observe first that any (bi)secant line to 
$C_H$ through a general point in $L$ 
meets $\Sigma_H$ in at least $k+k+1$ points counting multiplicities. 
Since $\Sigma_H$ has degree 
$a=k+1<2k+1$, it follows that such a line is contained in  $\Sigma_H$. If $C_H$ is not a plane curve,
it follows that through a general point of it there pass at least two lines of it (the line $L$ and one of 
the bisecant lines mentioned above). But it is well known that a ruled surface with such a property is necessary 
a quadric surface, which contradicts Lemma~\ref{lem:quadrica}. 
 
Therefore $F$ is a degenerate surface; actually, it must be a plane: in fact, 
if $P$ is a general point in $\gen{C_H}$, then through it there passes 
only one line $\ell_P$ of the congruence, and its focal point 
$\ell_P\cap D$ is contained in $C_H$. 


By Theorem~\textup{8} 
of \cite{DP2}, in which we classified all first order 
congruences in $\p^4$ whose fundamental locus is set-theoretically linear, 
we finish the proof. 
\end{proof}


\subsection{The Case in which $\length(\Lambda(b)\cap (F)_{\red})=2$}

Of these congruences we give a complete classification in the following 

\begin{thm}\label{thm:l2}
If $B$ is a first order congruence of lines given by a subfamily of the secant lines 
to the reduced locus $D:=(F)_{\red}$ of its pure fundamental locus $F$, then $D$ is 
a non-degenerate cubic scroll in $\p^4$ and the congruence is given as in Example~\textup{\ref{ex:2}}. 

\textit{Vice versa} a $B$ constructed in this way is a first order 
congruence such that $\length(\Lambda(b)\cap D)=2$. 
\end{thm}

\begin{proof}
Let us consider the surface $\Sigma_H$ of degree $a$ associated to a general 
hyperplane $H$; it contains the curve $C_H$ with some (algebraic) multiplicity, say $k'$. 

As in the proof of the preceding theorem, we fix $L\subset\Sigma_H$, and $l\in L^\vee_H$.  
Then $l\cap\Sigma_H=L\cup C$,  where $C$ is a curve of degree $a-1$, 
with two points of multiplicity $k'-1$ in $L\cap F$. As before,  
$C$ intersects $L$ in only another point, by Corollary~\ref{lem:curva}, and
Lemmas~\ref{lem:cono} 
and \ref{lem:quadrica} and $a=2k'$. 

In this case,  
$C_H$ must be a twisted cubic: 
in fact, we recall that the twisted cubic is the sole irreducible curve 
with only one apparent double point, and 
so if $C_H$ is not the twisted cubic, the secant lines to $C_H$ passing through a point 
$P\in L$ which are distinct from $L$ are contained in $\Sigma_H$, 
by degree reasons, 
and so $\Sigma_H$ must be a quadric, contradicting Lemma~\ref{lem:quadrica}. 
Therefore we have only a secant line to $C_H$ through $P$, $L$, 
and $C_H$ can only 
be a twisted cubic and $(F)_{\red}$ a rational normal scroll of degree 
three in $\p^4$. We recall that the planes of the 
conics in the cubic scrolls cover $\p^4$: 
therefore the lines of the congruence contained in one of these planes can be a 
congruence of order one and we are in the case of the assertion of the theorem, or 
it is a finite set. But the last case cannot occur: in fact, every secant line 
determines in a unique way  a conic in the scroll, so we would have that 
$\dim(B)=2$. 

$k'=\deg(\phi_{\abs{C_0},F}):=k$ since $k'$ is the degree of the cone of the lines in $B$ through 
a general point in the scroll. So, by what we said above, $k=k'$ and $a=2k$. 

\textit{Vice versa}, it is easy to see that a congruence constructed in such a 
way, has order one. 
\end{proof}

\begin{rmk}
We note that the case of the cone in the preceding theorem 
was not considered by G.~Marletta in \cite{M3}. 

Moreover, the only possible smooth surface in Theorem \ref{thm:l2} is the 
rational normal cubic scroll in $\p^4$, so we have done case \eqref{1:1} of 
Theorem \ref{thm:1er}.
\end{rmk}

\section{The Case of a Reducible and Non-reduced Fundamental Surface}
\label{sec:3}

Now we consider the case in which $F$ is reducible but not (generically) 
reduced;
we have only one possibility: 
the general line $\Lambda(b)$ intersects a component $F_1$ of $F$ 
in a fat point of length two, and so $F_1$ is non-reduced, and the other component 
$F_2$ (which is reduced) of $F$ in a simple point.

We recall we have denoted $C_{H,1}:=(F_1)_{\red}\cap H$ and 
$C_{H,2}:=F_2\cap H$, so $C_H=C_{H,1}\cup C_{H,2}$ and $k_i$ is the algebraic 
multiplicity of $C_{H,i}$ in $\Sigma_H$.

We start with 

\begin{lem}\label{lem:k1k2}
Either $D:=(F_1)_{\red}$ or $F_2$ is a plane. Moreover, the class of the 
congruence is  $a=k_1+k_2$. 
\end{lem}

\begin{proof}
As in the proof of the preceding two theorems, we fix $L\subset\Sigma_H$, and $l\in L^\vee_H$.  
Then $l\cap\Sigma_H=L\cup C$,  where $C$ is a curve of degree $a-1$, 
with a point of multiplicity $k_1-1$ 
and a point of multiplicity $k_2-1$ in $L\cap F$. 
$C$ intersects $L$ in only another point, by Corollary~\ref{lem:curva}, and 
Lemmas~\ref{lem:cono} and \ref{lem:quadrica}, and $a=k_1+k_2$. 
Reasoning as before, through a general point in $L$ there will not pass 
another joining line $C_{H,1}$ and  $C_{H,2}$, and so these lines generate a
first order congruence in $H$. Then, since all the first order congruences of 
lines in $\p^3$ are classified (see, for example, \cite{DP4} Theorem~0.1), we have that 
one of these curves is a line. 
\end{proof}

\begin{prop}\label{prop:nr}
Let $B$ be a first order congruence with non-reduced and reducible fundamental locus $F$, 
whose non-reduced component is $F_1$ and the other is $F_2$. 
Then we have the following possibilities:
\begin{enumerate}
\item $D:=(F_1)_{\red}$ is a plane, and we have the following cases:
\begin{enumerate}
\item\label{eq:prim} $F_2$ is a non-degenerate rational scroll of type $(1,d)$ and  
$D\cap F_2$ is a unisecant line in $F_2$ and the congruence is as in Example~\textup{\ref{ex:3}}; or 
\item\label{eq:second} $F_2$ is a non-degenerate rational surface with sectional genus 
zero (\ie either a rational scroll or the projected Veronese surface) and the  congruence is 
as in Example~\textup{\ref{ex:4}};  or 
\end{enumerate}
\item\label{eq:ter} $F_2$ is a plane, and $D=(F_1)_{\red}$ is a non-degenerate rational 
scroll of type $(1,c)$, 
$D\cap F_2$ is a unisecant curve in $D$ of degree $c$ and  the congruence is as in Example~\textup{\ref{ex:5}}.  
\end{enumerate}
\end{prop} 

\begin{proof}
As we observed in the remark in Section~\ref{sec:2b}, if we have a linear 
component in the pure fundamental locus, then 
the lines of $B$ contained in a general hyperplane containing this component  
give a first order congruence in $\p^3$. 
So, if $F_2$ is this plane,  
the congruence induced in an $f_2\in F_2^\vee$ 
has a non-reduced line as a fundamental locus, 
given by $F_1\cap f_2$ 
(see for example Theorem~0.1, case~(1b) of \cite{DP4}), from which we obtain 
case~\eqref{eq:ter}. 

If instead $D=(F_1)_{\red}$ is a plane, then the congruence induced in 
a $d\in D^\vee$ 
is a congruence with two fundamental curves, $L$ and $C$, $L$ is a line and 
if $\deg(C)=c$, then $\length(C\cap L)=c-1$ (Theorem~0.1, case~(2) 
of \cite{DP4}).  
If $L\subset F_2$, we have case~\eqref{eq:prim}, while if 
$L\subset (F_1)_{\red}$ we have case~\eqref{eq:second}. 

The calculations of $k_1$ and $k_2$ are immediate once one remembers that 
these numbers are the degrees of the cones of the lines of $B$ passing through 
a general point in, respectively, $F_1$ and $F_2$. The class follows from 
Lemma~\ref{lem:k1k2}. 
\end{proof}

From this, we easily obtain Theorem~\ref{thm:1er}, case~\eqref{1er}.

\providecommand{\bysame}{\leavevmode\hbox to3em{\hrulefill}\thinspace}
\providecommand{\MR}{\relax\ifhmode\unskip\space\fi MR }
\providecommand{\MRhref}[2]{%
  \href{http://www.ams.org/mathscinet-getitem?mr=#1}{#2}
}
\providecommand{\href}[2]{#2}


\begin{thebibliography}{{De }03b}

\bibitem[AF1]{AF}
S.~I. Agafonov and E.~V. Ferapontov, \emph{Systems of conservation laws from
  the point of view of the projective theory of congruences}, Izv. Ross. Akad.
  Nauk Ser. Mat. \textbf{60} (1996), no.~6, 3--30.

\bibitem[AF2]{AF1}
\bysame, \emph{Theory of congruences and systems of conservation laws}, J.
  Math. Sci. (New York) \textbf{94} (1999), no.~5, 1748--1794, Geometry, 4.

\bibitem[AF3]{AF2}
\bysame, \emph{Systems of conservation laws of Temple class, equations of
  associativity and linear congruences in $\mathbb{P}^4$}, Manuscripta Math.
  \textbf{106} (2001), no.~4, 461--488.

\bibitem[{De }99]{tesi}
P.~{De Poi}, \emph{On First Order Congruences of Lines}, Ph.D. thesis, SISSA-ISAS,
  October 1999.

\bibitem[{De }01]{DP2}
\bysame, \emph{On first order congruences of lines of $\p^4$ with a fundamental
  curve}, Manuscripta Math. \textbf{106} (2001), 101--116.

\bibitem[{De }03]{DP3}
\bysame, \emph{Threefolds in $\mathbb{P}^5$ with one apparent quadruple point},
  Comm. Algebra \textbf{31} (2003), no.~4, 1927--1947.

 \bibitem[{De }04]{DP4}
 {\bysame},
 \emph{Congruences of lines with one-dimensional focal locus},
 {Portugaliae Mathematica (N.S.)}  \textbf{61}
 ({2004}), no.~3, 329--338.

\bibitem[{De }05]{DPi}
\bysame, \emph{On first order congruences of lines in $\mathbb{P}^4$ with
  irreducible fundamental surface},
{Mathematische Nachrichten} \textbf{278} (2005), no.~4 363--378.

\bibitem[GH78]{GH}
P.~Griffiths and J.~Harris, \emph{Principles of Algebraic Geometry}, John Wiley
  \& Sons, 1978.

\bibitem[Har77]{H}
R.~Hartshorne, \emph{Algebraic Geometry}, Graduate Texts in Mathematics,
  no.~52, Springer-Verlag, 1977.

\bibitem[Kum66]{K}
E.~E. Kummer, \emph{\"{U}ber die algebraischen {S}trahlensysteme, insbesondere
  \"uber die der ersten und zweiten {O}rdnung}, Abh. K. Preuss. Akad. Wiss.
  Berlin (1866), 1--120, also in E. E. Kummer, \emph{Collected Papers},
  Springer Verlag, 1975.


\bibitem[Mar09a]{M3}
G.~Marletta, \emph{Sopra i complessi d'ordine uno dell'${S}_4$}, Atti Accad.
  Gioenia, Serie V, Catania \textbf{III} (1909), 1--15, Memoria II.

\bibitem[Mar09b]{M1}
\bysame, \emph{Sui complessi di rette del primo ordine dello spazio a quattro
  dimensioni}, Rend. Circ. Mat. Palermo \textbf{XXVIII} (1909), 353--399.

\bibitem[MP97]{MP}
E.~Mezzetti and D.~Portelli, \emph{A tour through some classical theorems on
  algebraic surfaces}, An. \c{S}t. Univ. Ovidius Contan\c{t}a \textbf{5}
  (1997), no.~2, 51--78.

\bibitem[O]{O}
A. Oblomkov, \emph{private communication}.

\bibitem[Ran86]{R}
Z.~Ran, \emph{Surfaces of order $1$ in {G}rassmannians}, J. Reine Angew. Math.
  \textbf{368} (1986), 119--126.

\bibitem[Seg88]{S}
C. Segre, \emph{Un'osservazione sui sistemi di rette degli spaz\^\i\ 
superiori}, {Rend. Circ. Mat. Palermo} \textbf{II} (1888), 148--149.



\bibitem[Zak93]{Zl}
F.~L. Zak, \emph{Tangents and Secants of Algebraic Varieties}, American
  Mathematical Society, Providence, RI, 1993, translated from the Russian
  manuscript by the author.

\end{thebibliography}
\end{document}